\documentclass{article}

\usepackage{arxiv}

\usepackage[utf8]{inputenc} 
\usepackage[T1]{fontenc}    
\usepackage{hyperref}       
\usepackage{url}            
\usepackage{booktabs}       
\usepackage{amsfonts}       
\usepackage{nicefrac}       
\usepackage{microtype}      
\usepackage{lipsum}         
\usepackage{graphicx}
\usepackage{natbib}
\usepackage{doi}

\usepackage{enumerate}
\usepackage{tikz}
\usepackage{float, upgreek}
\usepackage{bbm}
\usepackage{amsmath, amsfonts, amsthm, mathrsfs, amssymb}
\usepackage{appendix}
\usepackage{xcolor}
\usepackage{url}
\usepackage{subfig}
\usepackage{multirow} 
\usepackage{tweaklist}
\usepackage{cleveref}

\newtheorem{theorem}{Theorem}[section]

\newtheorem{lemma}[theorem]{Lemma}
\newtheorem{assumption}[theorem]{Assumption}
\newtheorem*{remark}{Remark}
\usepackage[normalem]{ulem}

\DeclareMathOperator*{\argmin}{\arg\!\min}
\DeclareMathOperator{\Tr}{Tr}
\DeclareMathOperator*{\argmax}{\arg\!\max}

\title{Regret Bounds and Experimental Design for Estimate-then-Optimize}

\author{ Samuel Tan\\
	School of Operations Research and Information Engineering\\
	Cornell University\\
	\And
Peter Frazier\\
	School of Operations Research and Information Engineering\\
	Cornell University\\ \\
}

\renewcommand{\headeright}{}

\renewcommand{\shorttitle}{Regret Bounds and Experimental Design for Estimate-then-optimize}

\hypersetup{
pdftitle={Regret Bounds and Experimental Design for Estimate-then-Optimize},
pdfsubject={math.oc, stat.ml},
pdfauthor={Samuel Tan, Peter Frazier},
pdfkeywords={Estimate-then-optimize, Experimental Design, Optimization Under Uncertainty},
}

\begin{document}
\maketitle

\begin{abstract}
  In practical applications, data is used to make decisions in two steps: estimation and optimization. First, a machine learning model estimates parameters for a structural model relating decisions to outcomes. Second, a decision is chosen to optimize the structural model's predicted outcome as if its parameters were correctly estimated. 
  Due to its flexibility and simple implementation,
  this ``estimate-then-optimize'' approach is often used for data-driven decision-making.
  Errors in the estimation step can lead estimate-then-optimize to sub-optimal decisions that result in regret, i.e., a difference in value between the decision made and the best decision available with knowledge of the structural model's parameters.
 We provide a novel bound on this regret for smooth and unconstrained optimization problems. Using this bound, in settings where estimated parameters are linear transformations of sub-Gaussian random vectors, we provide a general procedure for experimental design to minimize the regret resulting from estimate-then-optimize. We demonstrate our approach on simple examples and a pandemic control application.
\end{abstract}


\section{Introduction}
Consider a decision $x\!\in\!\mathbb{R}^d$ based on a model whose parameter vector $\theta$ is estimated from data. Examples include:
\begin{itemize}
\item Choosing nonpharmaceutical interventions to prevent the spread of a virus given a structural epidemiological model of viral spread, where the model's predictions depend on viral transmission parameters;
\item Choosing staffing levels in a hospital to ensure timely service and good health outcomes while controlling costs, given a structural queuing model for patient arrivals and the time required to deliver care;
\item Choosing rates at which to pump oil from several wells to maximize the amount of oil extracted, given a structural model for subsurface oil flow.
\end{itemize}

Perhaps the most common approach to solving this problem is to (1) estimate the parameter vector $\theta$ from data; and then (2) choose decision vector $x$ to maximize the value predicted by the structural model assuming this estimated parameter vector is correct. 
Mathematically, if our estimator is $\hat{\theta}$ and our model is $f(\cdot,\cdot)$, the resulting decision is $\hat{x} \in \argmax_x f(x,\hat{\theta})$.
This approach is sometimes called ``estimate-then-optimize'' \citep{farias2007revenue,hu2020fast}.

While more sophisticated methods exist for choosing a decision $\hat{x}$ based on data, such as Bayesian methods \citep{wu2018bayesian} and robust optimization \citep{bertsimas2018data}, they can be hard to implement. One barrier is 
the computation required by the more complex optimization problems they solve (e.g., $\max_x \min_{\hat{\theta}} f(x,\hat{\theta})$ for robust optimization).
Another barrier 
is the need to integrate engineering and organizational structures that previously separated estimation from optimization. For example, a company may have a team of machine learning experts focused on estimation and a separate team trained in and focused on optimization, each with their own engineering systems.
Thus, estimate-then-optimize is widely used in practice and it is valuable to understand and improve its performance.

In estimate-then-optimize, errors in the statistical estimation can lead to regret: the value of the decision made $f(\hat{x},\theta)$ may be smaller than that of the best decision possible using the true parameters, $\max_x f(x,\theta)$. While this is understood intuitively to be important by practioners, past literature has not, to our knowledge, provided formal bounds for this regret nor considered how the experimental design used to gather data for estimation impacts regret.

Our contribution is a novel bound on the regret resulting from an estimate-then-optimize approach for smooth $f$ and unconstrained $x$. 
We further show how this approach can inform the experimental design used to collect data for estimating $\theta$. In particular, our results can help a user decide how much effort to allocate to parts of their parameters space so as to minimize regret. In a setting where the parameter estimate is known to be a linear transformation of a sub-Gaussian random vector, we provide bounds that depend on its covariance; this covariance term can consequently be used to minimize this asymptotic regret bound. We provide three examples: one in which the parameters are directly assumed to be normal, a pricing example using a logistic conversion model, and a more tangible real-world example of controlling a pandemic.

We envision this work will be of substantial interest to those who currently rely on the estimate-then-optimize paradigm to make decisions from data and are prevented by the barriers described above from moving away from estimate-then-optimize. For instance, a marketplace platform may use our result to select not only prices to use for an experiment, but also the proportion of participants to assign to each price so as to minimize the regret in finding the optimal price.





The rest of the paper is structured as follows. Section \ref{sec:rel} overviews related work. Section \ref{sec:tech} contains the main technical results and the prerequisite assumptions. Section \ref{sec:norm}  expands on said technical results in cases when the parameter estimates are assumed to be linear transformations of sub-Gaussian random vectors. Section \ref{sec:app} connects our work to experimental design. Section \ref{sec:exam} demonstrates our method on three examples. Finally, section \ref{sec:future} concludes and mentions some future directions of this work.

\section{Related Work}
\label{sec:rel}

The general two-stage procedure of estimating parameters in an optimization problem and then solving using those predicted parameters has appeared in the literature in various forms.

From a statistical viewpoint, the framework presented in this paper can be viewed as optimizing some function in the presence of a model parameter (sometimes called a nuisance parameter), $\theta$, which must first be estimated. Using this perspective, \cite{foster2019orthogonal} derive fast rates of convergence 
when a condition called Neyman-orthogonality is satisfied, 
where the pathwise derivatives of $f$ with respect to $\theta$ and $x$ at their true values
vanishes. However, this cross-derivative term is exactly the term which we consider in this work as (1) it is generally not equal to 0 except in very specific applications and (2) it governs the slower rate for the risk.

There has also been work done in similar settings in which one wishes to perform estimate-then-optimize in a ``smarter'' way. In \citep{elmachtoub2022smart}, the authors specifically tackle a contextual stochastic optimization problem, wherein the objective is linear and of the form $\min_{x} \mathbb{E}_{c(w) \sim F_w}[c(w)^\top x | w]. $ Here the authors call this particular form ``predict-then-optimize''.
In particular, here the user is presented with a context $w$ drawn from some distribution, and the parameters $c(w)$ of the objective function depend on $w$. The standard estimate-then-optimize approach would be to run some off-the-shelf regression method to construct a prediction $\hat{c}$ for $\mathbb{E}[c(w)|w]$ based on past data of the form $\{(c_i, w_i)\}_{i=1}^n$, and optimize the resulting program $\min_{x} \hat{c}^\top x$. The proposed solution is to aim the prediction of $\mathbb{E}[c(w)|w]$ towards minimizing the objective value in the downstream optimization problem, dubbed the smart predict-then-optimize (SPO) loss. As minimizing this directly is intractable, the authors propose minimizing a convex surrogate. Later work, such as \cite{el2019generalization}, analyze generalization and risk bounds for this particular approach.

We emphasize that this ``predict-then-optimize'' framework focuses on the loss for individual contexts, i.e. $\mathbb{E}[\theta|w] \cdot x(w)$ for a fixed $w$. However, if one instead focuses on the average performance across $w_i$ chosen from the training set, then the overall loss would be of the form $\sum_{i=1}^n \mathbb{E}[c(w_i)|w_i] \cdot x(w_i)$. We may encode this in our framework via setting $\theta:= \{ \mathbb{E}[c(w)|w]: w \}$ and $x := \{ x(w) : w\}$, and $f(x,\theta) := \sum_{i=1}^n \theta(w_i) x(w_i)$.

\cite{ito2018unbiased} also consider the same estimate-then-optimize framework, wherein one optimizes $f(x,\hat{\theta})$ using an estimate $\hat{\theta}$ for the true $\theta^*$. However, their work does not analyze the value of $f(\hat{x}, \theta^*)$, but rather establishes that the estimate $f(\hat{x}, \hat{\theta})$ is generally biased for $f(x^*, \hat{\theta})$. They additionally describe a bias-correction procedure under $f$ which is affine in $\theta$, but do not offer asymptotic analysis of our quantity of interest, $f(\hat{x}, \theta^*) - f(x^*, \theta^*)$.

\cite{wilder2019melding} also consider the estimate-then-optimize framework, though in their case the optimization problem is combinatorial, i.e. the decision $x$ are discrete. Their approach is to (1) relax the combinatorial problem to a continuous one, (2) derive gradients of the objective with respect to the predictions, (3) use the closed-form gradients in a gradient-based optimization algorithm, and (4) round the resulting solution to the nearest integer. They also find that $\frac{\partial^2 f}{\partial x \partial \theta}$ plays a fundamental role in their gradient derivations, and demonstrate their approach in detail on linear programming and submodular maximization problems. However, they only offer empirical evidence that their approach performs well, and obviously apply only to combinatorial problems.

Our work is also closely related to empirical risk minimization (ERM) \citep{vapnik1999nature} and sample average approximation (SAA) \citep{kim2015guide}, though still different in important ways. In our setting, toward the goal of minimizing $f(x,\theta)$ over $x$, we consider the approximate minimizer $\hat{x} \in \argmax_{x} f(x,\hat{\theta})$ for some estimator $\hat{\theta}$ of $\theta$. In ERM and SAA, toward the goal of minimizing $\mathbb{E}_\psi g(x,\psi)$ over $x$, one considers the approximate minimizer $\hat{x} \in \argmax_x \frac{1}{n} \sum_{i=1}^n g(x,\psi_i)$ where $\psi_i$ are i.i.d. Except in trivial cases (such as those where $f,g$  do not depend on  $\theta$), there is no way to construct an $f$ and $g$ so that the quantities match. Our setting is therefore  different from ERM/SAA except in the special case where $f$ is linear in $\theta$ for all $x$.

The idea of using sensitivity of output to parameters to guide data collection has also been explored, notably in the simulation community, e.g. \cite{song2015quickly} and \cite{freimer2002collecting} give two approaches for deciding on which parameters to collect more data in order to reduce the impact of input uncertainty on simulation output. However, these approaches do not consider the setting in which there is a decision to be optimized, whereas our approach is concerned with targeting parameters to make better decisions.

Finally, experimental design is a mature field within statistics which has explored the notion of \textit{optimal} design for an experiment, where there are several notions of optimality \citep{pukelsheim2006optimal}. In the particular case of linear models, in which the experimenter has control over the proportion of samples to allocate to a set of linearly independent vectors which comprise the design matrix, various optimality criteria have been studied in the literature. A commonly used definition of optimality is $A$-optimality, which reduces to minimizing the average variance of the regression coefficients. Similarly, $C$-optimality instead considers minimizing the variance of a linear combination of the regression coefficients. When the decision variable is one-dimensional, our regret bound induces a $C$-optimal design problem. More generally, our regret bound induces an $A$-optimal design for a particular parameter system $K\theta$.

\section{General Regret Bound for Estimate-then-Optimize} 
\label{sec:tech}
This section defines the estimate-then-optimize framework more formally and then states our general bound on regret and assumptions required by this bound.

In the estimate-then-optimize framework we are concerned with
optimizing some smooth function $f(x,\theta)$ 
where $\theta$ is a model parameter
and $x \in \mathbb{X}$. Here $\mathbb{X}$ can equal $\mathbb{R}^d$ or a subset thereof.
The goal is to solve $\argmin_x f(x,\theta^*)$,
where $\theta^* \in \Theta$ is unknown and must be estimated.
Call $x^* = \argmin_x f(x,\theta^*)$. Note that here if $\mathbb{X} \subsetneq \mathbb{R}^d $, Assumption 3.1 must be satisfied, i.e. $x^* \in \text{int } \mathbb{X}$.

Below, we denote $\mathcal{S}(z, z') := \{tz + (1-t)z'| t\in [0,1] \}$, i.e. the line segment between $z, z'$. Also $\| \cdot \|$ (without a subscript) always denotes the Euclidean norm. Finally, for higher-than-second-order directional derivatives (as they appear in Taylor's theorem), we denote $D^k_s f(s)[p_1, \ldots, p_k]$ to be 
\[ \frac{\partial^k}{\partial t_1 \ldots \partial t_k}f(s+ t_1 p_1 + \cdots + t_kp_k)\Bigg|_{t_1 = \cdots = t_k = 0} .\]

We consider an algorithm that first learns $\hat{\theta}$ from a sample $S$ 
(say a sample generated under a distribution with mean $\theta^*$)
and then learns $\hat{x}$ by minimizing $f(x,\hat{\theta})$, where the second 
minimization step is exact.


\begin{assumption}[First Order Optimality]
	The minimizer $x^*$ satisfies 
	$ \frac{\partial f}{\partial x} \Big|_{x = x^*, \theta = \theta^*} =0 .$ 
\end{assumption}

\begin{assumption}[Strong Convexity]
The function $f$ is strongly convex with respect to $x$ in a neighorhood $B_\theta$ around $x^*(\theta) := \argmin_x f(x,\theta)$ for each $\theta$, i.e. a universal $\rho$ exists such that 

\[ (x- x^*)^\top \frac{\partial^2 f}{\partial x^2} \Bigg|_{x=\overline{x}, \theta = \theta'} (x-x^*) \geq \rho \|x-x^*\|^2 \]
for all $x \in B_{\theta'}, \overline{x} \in \mathcal{S}(x, x^*)$, for all $\theta'$.
\end{assumption}

\begin{assumption}[Higher Order Smoothness]
	There exist constants $\beta_1$ and $\beta_2$ such that the following
	derivative bounds hold:
	
	\begin{itemize}
		\item For all $x \in \mathcal{X}$, $\overline{x} \in \mathcal{S}(x, x^*)$
			\[  (x-x^*)^\top \frac{\partial^2 f}{\partial x^2} \Bigg|_{x=\overline{x}, \theta = \theta^*} (x-x^*) \leq \beta_1 \|x - x^* \|^2,\] 
		\item For all $\theta \in \Theta$, $\overline{\theta} \in \mathcal{S}(\theta, \theta^*)$, $x \in \mathcal{X}$
			\begin{align*} \Big | D^2_\theta D_x f(x^*, \overline{\theta}) [x - x^*, \theta - \theta^*, \theta - \theta^*] \Big | \\
			\leq \beta_2 \| x - x^*\| \|\theta  - \theta^* \|^2 .\end{align*}
		\end{itemize}
\end{assumption}

The higher order assumptions also appear in other literature, e.g. \cite{foster2019orthogonal}, where it is shown that they are easily satisfied. We then have the following theorem.

\begin{theorem}
	Under the above assumptions and when $x^* \in B_{\hat{\theta}}$,
	the optimality gap of $\hat{x}$ under the true $\theta^*$, i.e. $f(\hat{x}, \theta^*) - f(x^*, \theta^*) $, is upper bounded by $\frac{4\beta_1}{\rho^2} \Bigg( \|\frac{\partial^2 f}{\partial x \partial \theta}\Bigg|_{x = x^*, \theta = \theta^*} (\hat{\theta} - \theta^*) \| ^2 +  \frac{\beta_2^2}{4} \|\hat{\theta} - \theta^* \|^4 \Bigg).$
	\label{thm:main}
\end{theorem}

\section{Regret Bound for Linear Transformations of sub-Gaussian Random Vectors $\hat{\theta}$}
\label{sec:norm}

We now consider a particular setting in which we place a distributional assumption on our estimate $\hat{\theta}$. In particular, suppose we have some dataset of size $n$ which we use to construct $\hat{\theta}$. Furthermore, assume this estimate is known to be a linear transformation of a sub-Gaussian random vector which has parameter 1 WLOG. This encompasses a wide class of random variables. In particular, normal random variables are a special case when taking the sub-Gaussian random vector to be an isotropic multivariate Gaussian with identity covariance and the linear transformation to be the square root of the covariance matrix. In the asymptotic regime, asymptotically normal estimators are therefore also included, e.g. from the central limit theorem or asymptotic theory governing maximum likelihood estimators. This is summarized in the following assumption.

\begin{assumption} 
\label{as:norm}
We have an algorithm which, from a dataset $\mathcal{S}$ of size $n$, produces $\hat{\theta}$ such that $\hat{\theta}-\theta =  \frac{\Sigma^{1/2}}{\sqrt{n}} X$ for positive definite $\Sigma$ and $X \sim$ \text{subG(1)}.
\end{assumption}

\begin{assumption}
\label{as:open}
The set $\Theta^* := \{\theta :  x^*(\theta^*) \in B(\theta)\}$ is open.
\end{assumption}

Also for notational convenience we define $D:= \frac{\partial^2 f}{\partial x \partial \theta} \Big|_{x = x^*, \theta = \theta^*} .$ 
We then have the following lemma for converging to the strongly convex neighborhood.

\begin{lemma}
There exist constants $\alpha, \beta >0$ such that for $n$ sufficiently large,
$x^* \in B_{\hat{\theta}}$ with probability at least $1-\alpha \exp(-\beta n)$.
\end{lemma}





This lemma then yields the following theorem.
\begin{theorem}
\label{thm:hpb}
With probability at least $1- \alpha \exp(-\beta n) - \frac{2}{n}$
\begin{align*} &f(\hat{x}, \theta^*) - f(x^*, \theta^*)
\le  \frac{4\beta_1}{\rho^2} \Bigg( \Tr  \frac{D\Sigma D^\top}{n} \\
& + 2 \sqrt{\Tr \left[ \left(\frac{D\Sigma D^\top}{n}\right)^2\right] \log n } + 2 \|\frac{D\Sigma D^\top}{n}\|_2 \log n \Bigg)\\
& + o(\Tr \frac{D\Sigma D^\top}{n} \log n)
\end{align*}
\end{theorem}

\section{Application to Experimental Design}
\label{sec:app}

We now imagine we have a class of experiments which informs us about the uncertain parameter $\theta$. Each produces a different estimator and corresponding distribution, which is assumed to satisfy the distributional assumptions of Section \ref{sec:norm}. Our goal is to choose the experiment so as to minimize the upper bound derived previously in Theorem \ref{thm:hpb}. 

In particular, the theorem tells us that minimizing
the trace of the variance term $\frac{D \Sigma D^\top}{n}$ minimizes the upper bound on the regret. In the context of linear models, this is exactly equivalent to the $A$-optimality criterion in experimental design \citep{pukelsheim2006optimal} when using $D$ as a coefficient matrix. \footnote{Using other scalar values related to the matrix $D\Sigma D^\top$ would correspond to other criteria. For instance, minimizing the determinant (or equivalently log determinant), would correspond to $D$-optimality, or minimizing the 2-norm would correspond to $E$-optimality.}

Therefore the problem is to minimize $ \Tr \frac{D\Sigma D^\top}{n}$, where the decision variable is the experimental design which induces the form of $\hat{\theta}.$




In the discussion and results above, we assume $D$ is known to us. However in practice, $D$ is unknown --- it depends on the unknown value $\theta^*$. We overcome this with a Bayesian approach: we place a prior $\pi$ on $\theta^*$, and draw samples from the prior $\theta$, which we then use to calculate $x^*(\theta)$. This then gives us access to $D' := \frac{\partial^2 f}{\partial x \partial \theta} \Big|_{x^*(\theta), \theta}. $ We then solve the following problem

\begin{equation} \min_{\hat{\theta}} \mathbb{E}_{\theta \sim \pi} \left[  \Tr \frac{D'\Sigma D'^\top}{n} \right],\label{eq:bayes} \end{equation}
where in practice we replace the true expectation with an empirical average over draws from the prior.

\section{Examples}

This section provides three examples of how our method can be used to inform experimental design, as well as corresponding numerical experiments demonstrating the value of our method.

\label{sec:exam}

\subsection{Experimental Design Under Normal Model Parameters from a Linear Model}
\label{sec:example-linear}

We consider a particular setting for experimental design to which 
our method directly applies. Suppose
we can get samples $\theta_{ij} \sim N(\theta^*_i, \sigma_i^2),$
where $\theta_{ij}$ are i.i.d. for each $i$ and $\theta_i$ are independent between $i=1,\dots, n$.

First, for $i \in \{1, \ldots, n\}$, define $N_i$ as the number of times $\theta_{i,j}$ was sampled. Then we define $\hat{\theta}$ component-wise via
\[\hat{\theta}_i = \frac{\sum_{j=1}^{N_i}\theta_{i,j}}{N_i}. \] 
So in this case $\hat{\theta} - \theta \sim N(0, \hat{\Sigma}),$where $\hat{\Sigma}$ is diagonal with $\hat{\Sigma}_{ii} = \frac{\sigma_i^2}{N_i}.$

We see that minimizing the norm of the LHS is equivalent to minimizing the variance $\mathbf{d}^\top \hat{\Sigma} \mathbf{d}$ (for instance, the squared norm of the LHS is a $\chi^2$ distribution).

Now assume that $x \in \mathbb{R}$. Since $D$ in this case is just a vector, we then denote $ \mathbf{d}  = \frac{\partial^2 f}{\partial x \partial \theta}\Big|_{x = x^*, \theta = \theta^*}.$

Now referring back to Theorem \ref{thm:hpb}, we see that the quantity of interest is $\mathbf{d}^\top \Sigma \mathbf{d}$ (since the trace and $\|\cdot \|_2$ in this case correspond to the same scalar). So the objective is to minimize $\mathbf{d}^\top \Sigma \mathbf{d}$, subject to $\Sigma$ belonging to the class of allowed by our 
allocations $N_i$.

At this point, we point out that this exactly coincides with $C$-optimal experimental design. In particular, this is the special case that our linearly independent vectors for samples are the unit vectors in $\mathbb{R}^n$ so that $X = I$. Then standard results in $C$-optimal experimental design \citep{pukelsheim2006optimal} yield that, given a mass of unit $1$ to allocate across the $n$ vectors, 
\begin{equation} n_i = \frac{d_i \sigma_i}{ \sum_{j=1}^n d_j \sigma_j }. \label{eq:closed} \end{equation}
In practice, if one had say $m$ samples to allocate across the $n$ components, we would assign $mn_i$ to component $i$, rounding when necessary.

We now provide a simple problem to demonstrate our method on this setting. Consider $f(x,\theta) = \frac{\theta_1}{2} x^2 + \theta_0 x,$ which satisfies the assumptions of the theorem when $\theta_1 > 0$. Suppose that we can sample $\theta_i \sim N(\theta_i^*, \sigma_i^2)$. 
Here $D = \left(1 , -\frac{\theta^*_0}{\theta^*_1} \right)$, and minimizing the bound corresponds to minimizing $\frac{\sigma_0^2}{n} + \frac{\sigma_1^2}{N-n}\left(\frac{\theta_0^*}{\theta_1^*}\right)^2,$ where $N$ is the total number of samples to allocate. Note that we do not have a constant term, as any samples used to estimate the constant term are wasted; we care about $f(\hat{x}, \theta^*)$, and the constant term in $f(x,\theta^*)$ does not interact with the decision variable $x$.

\begin{figure}
  \centering
  \subfloat[][Regret vs. \# samples measuring the slope\label{fig:quad}]{%
  \includegraphics[width = 0.49\textwidth]{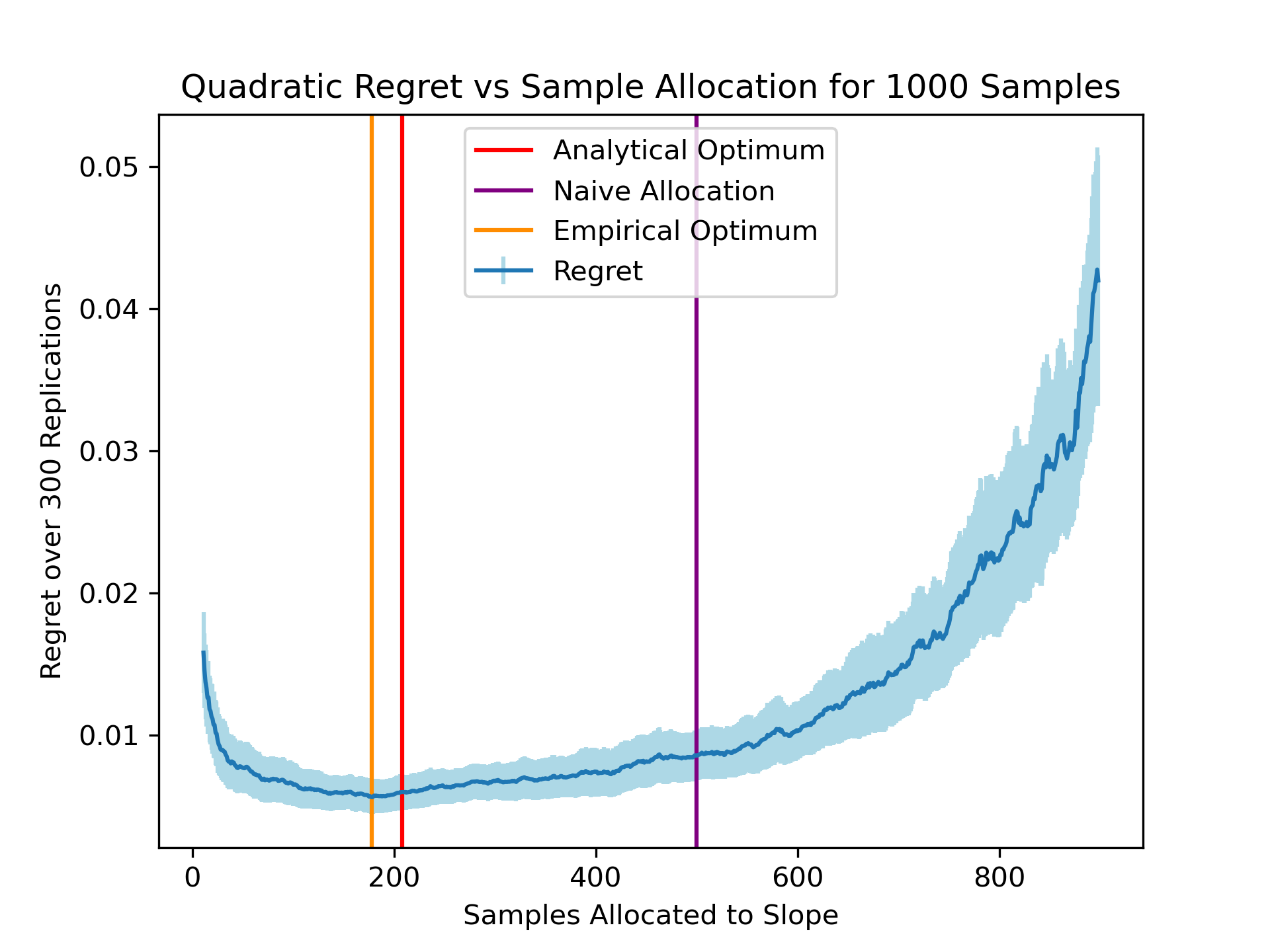}
  } \hfill
  \subfloat[Estimated objective function and decision \label{fig:effect}]{%
  \includegraphics[width = 0.49\textwidth]{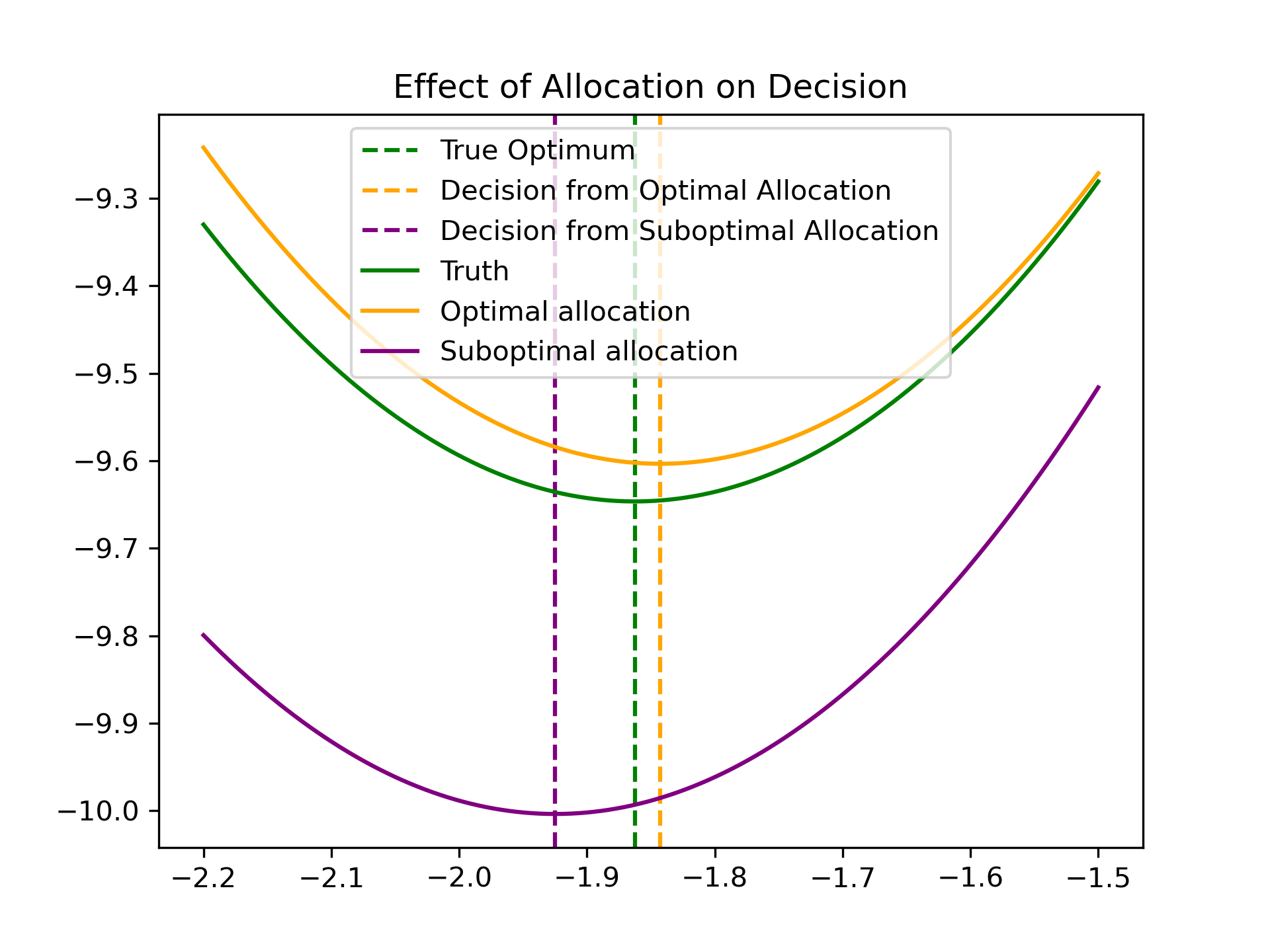}
  }
  \caption{Simple quadratic example illustrating the general setting in Section~\ref{sec:example-linear}. (a) A comparison of our method for collecting samples of $\theta_0$ and $\theta_1$. Our method corresponds to the analytical optimum. The empirical minimum regret (minimizer of the blue regret curve) is also shown, while the uniform allocation can be seen to achieve significantly higher regret.
  (b) The estimated quadratic function when using different allocations. The resulting decision from using the optimal allocation yields a solution closer to the optimizer of the true function.}
\end{figure}

For our experiment, we use the prior $\theta^* \sim N\Big( \begin{pmatrix}10 \\ 5 \end{pmatrix}, I\Big)$ and $\Sigma = \text{diag} (1,3)$.
We repeated the method for 300 replications, where 100 draws from the prior are used for both optimization of the bound and evaluation of the regret. Figure \ref{fig:quad} 
shows the regret under the prior achieved by our method versus 
a uniform allocation which splits samples across $\theta_0$ and $\theta_1$
evenly.

Figure \ref{fig:effect} compares the estimated function $f(x, \hat{\theta})$
against $f(x, \theta^*)$, showing that the resulting estimate $\hat{x}$ 
yields a better decision for $f(x, \theta^*)$ when using the optimal allocation.

\subsection{Pricing Under Logistic Conversion}
\label{sec:examp-price}
In this example, we consider pricing a product based on a data from a pricing experiment. We observe a sequence of customers arriving to a website or mobile app who are each presented with a price $x$. For each customer, we observe whether they purchase (``convert'') or not, with the probability of purchase being governed 
by a logistic function of the price $x$ with parameter $\theta$, i.e.

\[c(x,\theta) = \frac{1}{1  + \exp(\theta_0 + \theta_1x)} \]
The objective function is then revenue, 
negated to yield a minimization problem: $f(x,\theta) = -xc(x,\theta).$

Here we assume $\theta_1 > 1$ so that higher prices have lower purchase probability. We suppose that estimate-then-optimize is implemented by estimating $\theta$ via maximum likelihood estimation. From standard results for the maximum likelihood estimate for logistic regression \citep{hastie2009elements},
we know that asymptotically
$\hat{\theta} \sim N(\theta, (X^T W X)^{-1}),$
where $X$ is a $n \times 2$ matrix with first column
all ones, and second column corresponding to prices $x_i$,
and $W$ is a diagonal $n \times n$ matrix with $c(x_i, \theta)(1-c(x_i,\theta)$
on the diagonal.


In this case, $\mathbf{d}$ is given by
\begin{align*} \frac{\partial^2 f}{\partial x \partial \theta_0} =& c(x,\theta)(1-c(x,\theta))\\
&+ \theta_1 x \left(  c(x,\theta)(1-c(x,\theta) \left(2c(x,\theta) - 1 \right) \right) \\
\frac{\partial^2 f}{\partial x \partial \theta_1} =& 2x c(x,\theta)(1 - c(x,\theta))\\
&+ x^2 \theta_1 c(x,\theta)(1 - c(x,\theta)) \left(2c(x,\theta) - 1 \right)
\end{align*}

It is easy to show that $f$ satisfies the conditions of our theorem. For completeness, the proofs are provided in the supplementary material.

We consider a user who is designing an experiment to support selection of the best price. The user has access to historical data on prices and conversion, which she can use to construct a prior for $\theta$. Using this data, she designs an experiment using $n$ samples constrained to use prices from among $m$ distinct prices. The user then uses maximum likelihood to estimate $\theta$, and solving the corresponding revenue maximization problem. This yields the fixed price $\hat{x}$ which she holds in place moving forward. The revenue of the selected price $\hat{x}$, under the true logistic model $\theta^*$, is the quantity we wish to be as large as possible. The regret is the difference between this quantity and the revenue at the optimal price $x^*$ under the true logistic model.

We study this experimental design problem numerically via simulations.
We simulate the procedure described above, where one must decide how to allocate samples
(opportunities to sell the product to a customer) across prices to yield the selected price $\hat{x}$ with the highest expected revenue, considering $n=100$ samples and $m=10$ prices.

We evaluate our method and a baseline method under a collection of multivariate normal priors for $\theta^*$ with covariance $0.01I$.
The means of the priors for $\theta^*$ were chosen 
so that the $x$ value for which $c(x,\theta^*) = 0.5$
is at prices $1,2,\ldots, 8$ (avoiding boundaries). 100 draws from the prior were used for optimization and evaluation of regret.

Due to the computational cost of considering 
all partitions of 100 allocations among the 10 prices, we instead 
use random search to optimize the regret bound. In particular, we choose 1000 random allocations 
uniformly at random and return the allocation
minimizing $\frac{\mathbf{d}^\top (X^\top W X)^{-1} \mathbf{d}}{n}$.

We then simulate 
100 customer purchase decisions (conversions) using prices distributed among the 10 prices $0, \ldots, 9$ according to the best allocation found as well as according to the baseline uniform allocation, again using the logistic model and the held out ground truth $\theta^*$. We again use maximum likelihood to derive estimates $\hat{\theta}_{\text{optimized}}$ and $\hat{\theta}_{\text{uniform}}$, and their corresponding revenue maximizers 
$\hat{x}_{\text{optimized}}, \hat{x}_{\text{uniform}}$.

Table \ref{tab:res} gives 95\% confidence intervals 
for the resulting expected regrets from these two approximate revenue maximizers, 
$\hat{x}_{\text{optimized}}$ and  $\hat{x}_{\text{uniform}}$.
We see that our optimized approach provides consistently lower regret than the baseline uniform allocation. We also note that in our experiments, we observed the scale of the regret scales as $O(\frac{1}{n})$: see the supplementary material for more experiments for values of $n$.

\begin{table}\centering
  \begin{tabular}{cccc}\toprule
  &Optimized Allocation & Uniform Allocation \\
  $\theta = [-1,1]$ & $1.0 \pm 0.3$ & $1.4 \pm 0.4$ \\
  $\theta = [-2,1]$ & $0.8 \pm 0.3$ & $1.5 \pm 0.4$\\
  $\theta = [-3,1]$ & $1.2 \pm 0.3 $ & $1.5 \pm 0.4$ \\
  $\theta = [-4,1]$ & $1.3 \pm 0.3 $ & $1.8 \pm 0.5$ \\
  $\theta = [-5,1]$ & $2.4 \pm 2.3 $ & $2.5 \pm 0.6$ \\
  $\theta = [-6,1]$ & $2.1 \pm 0.6 $ & $3.3 \pm 1.0$ \\
  $\theta = [-7,1]$ & $2.8 \pm 0.8$ & $3.2 \pm 0.9$ \\
  $\theta = [-8,1]$ & $2.9 \pm 0.8$ & $6.6 \pm 2.8$\\
  \bottomrule
  \end{tabular}
  \caption{Comparison of regret under optimized and baseline uniform allocations in the pricing example of Section \ref{sec:examp-price}. The quantities are 95\% confidence intervals, where the units are in $10^{-2}$. 
  \label{tab:res}}
  \end{table}

\begin{figure}
  \centering
  \subfloat[][Revenue function and optimal allocation.\label{fig:price}]{%
  \includegraphics[width = 0.49\textwidth]{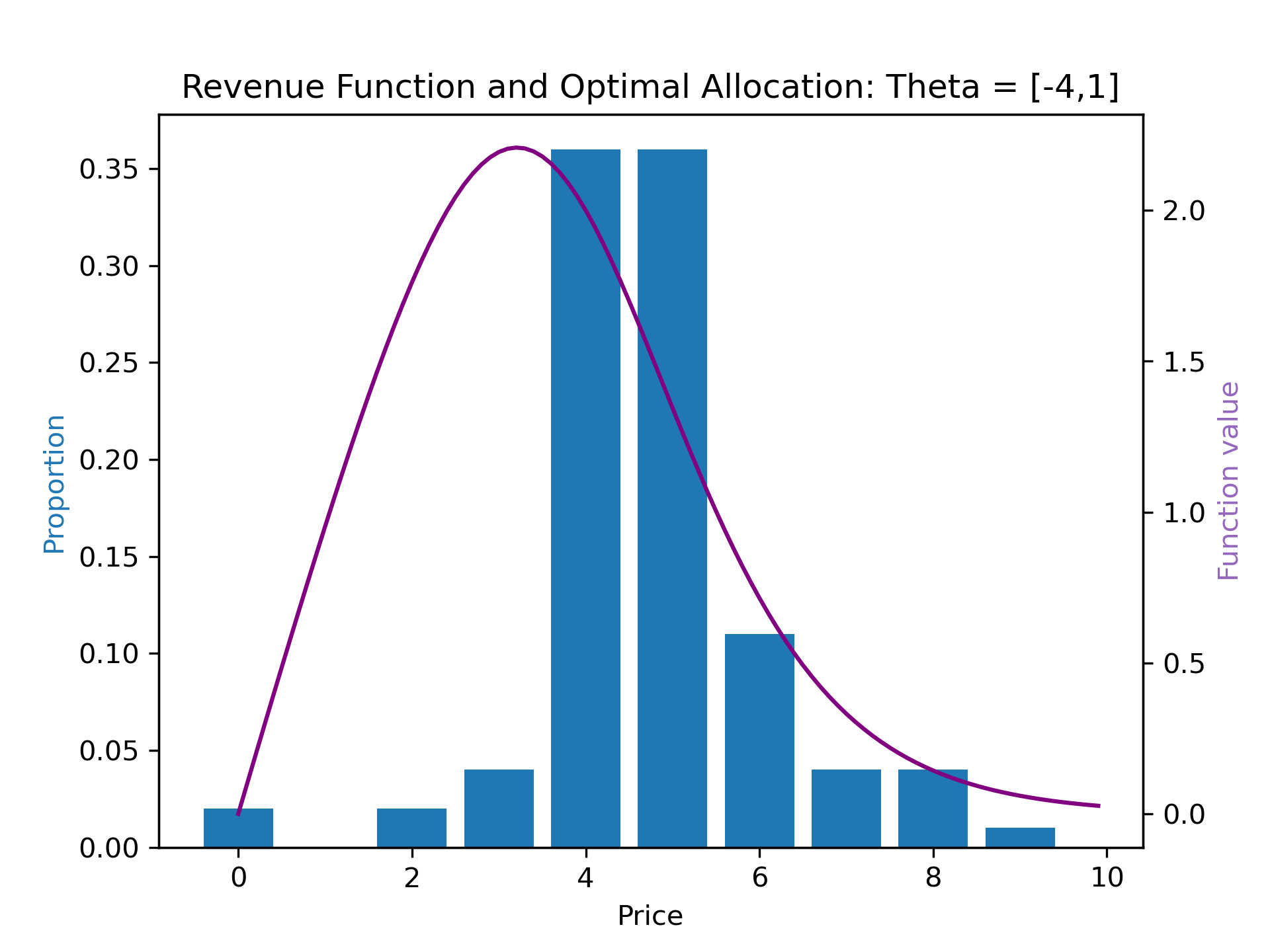}
  }\hfill
  \subfloat[Estimated objective function and decision \label{fig:price_eff}]{%
  \includegraphics[width = 0.49\textwidth]{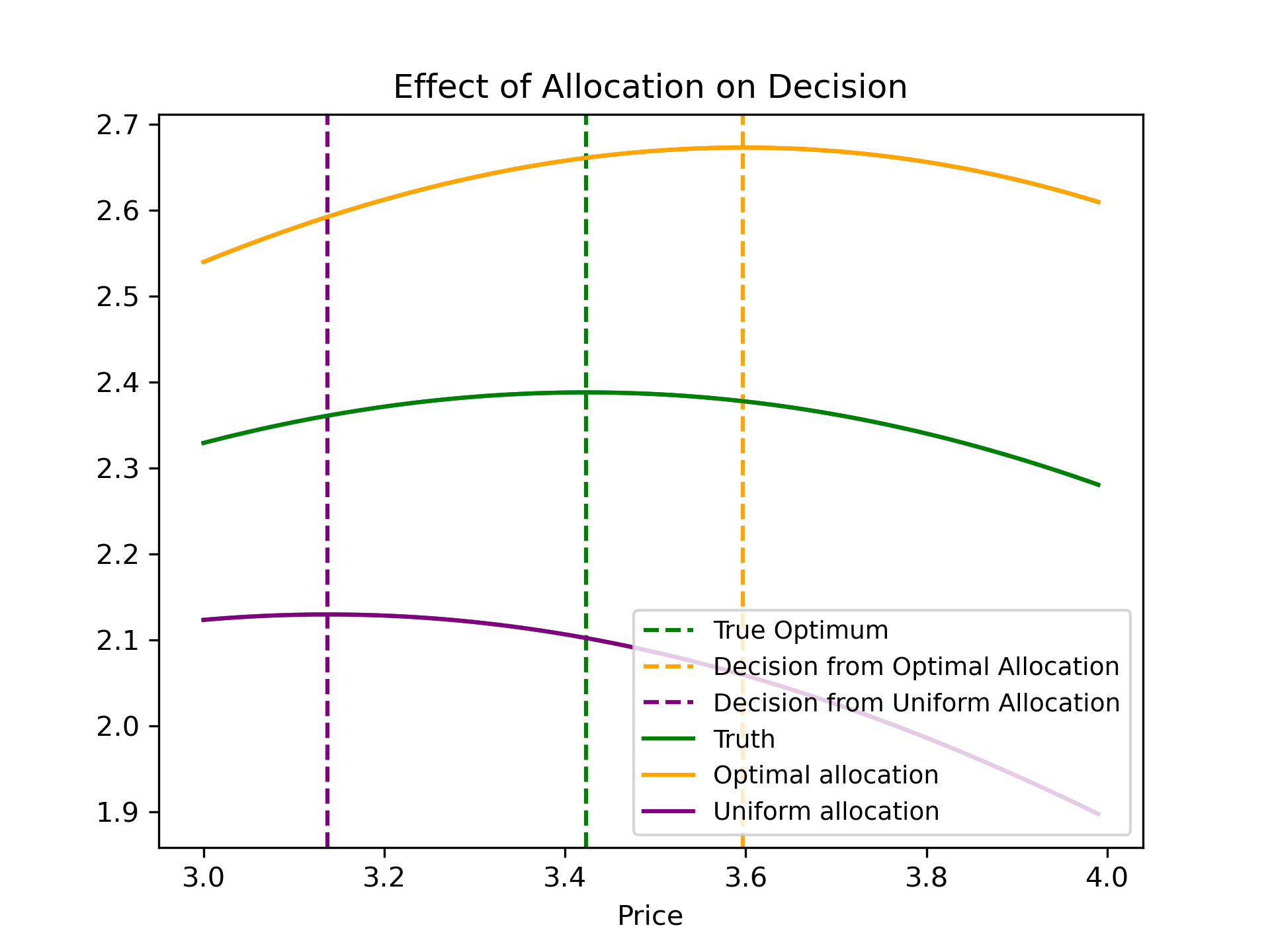}
  }
  \caption{Pricing example from Section \ref{sec:examp-price}. (a) True revenue function (purple line) and optimal experimental sample allocation (blue bars).
  (b) Estimated revenue function under three different experiment sample allocations, using a zoomed-in $x$-axis. The decision supported by the optimal allocation yields a solution closer to the optimizer of the true revenue.}
\end{figure}

  Figure \ref{fig:price} also shows an example of the optimal allocation
  found by our method. We can observe that, compared to a uniform sampling strategy, more samples
  are allocated near the peak of the function, yielding more information about the maximum. In particular,
  more samples are allocated slightly to the right where the conversion hovers around $\frac{1}{2}$ (at $x=4$), suggesting 
  that the change in curvature there also yields more information than 
  just the peak.
  
Finally, Figure \ref{fig:price_eff} is the analogous figure of Figure \ref{fig:effect} for this pricing problem.
Here it shows, again for $\theta^* = [-4,1]$, the estimated function $f(x,\hat{\theta})$
for $\hat{\theta}$ coming from different allocations. We again see that the optimal allocation
yields an estimated function whose optimizer is closer to the optimizer of the true function.

\subsection{Application to Pandemic Control}
\label{sec:pandemic}

We now apply our framework to controlling COVID-19 through regular testing as a university or city reopens, subject to uncertainty about contact among social groups.

We model viral spread via a multi-group SIR model (see, e.g., \citealt{guo2006global}), 
where people belong to one of $n$ social groups, indexed by $k$, and to one of three disease states.
A combination of social group and disease state is called a ``compartment''.
The three disease states considered are: Susceptible (S), someone who has never been infected; Infected (I), someone who is infected and can infect others; and Removed (R), someone who was previously infected, cannot infect others, and cannot be themselves infected. 
A person can be Removed either because they have recovered from their infection and developed immunity or because they were discovered through testing and isolated away from others.
$S_k$ indicates the number of susceptible people from social group $k$, and similarly for $I_k$ and $R_k$.

In this model, the following differential equations describe how people move between disease states:
\begin{align}
    \frac{dS_k}{dt} &= -\frac{\sum_{j=1}^n \beta_{kj} S_k I_j}{N_k}
    \label{eq:SIR-1}\\
    \frac{dI_k}{dt} &= \frac{\sum_{j=1}^n \beta_{kj} S_k I_j}{N_k} - (\gamma+t_k) I_k
    \label{eq:SIR-2}\\
    \frac{dR_k}{dt} &= (\gamma+t_k) I_k,
    \label{eq:SIR-3}
\end{align}
where $N_k$ is the size of the $k$th group, $\beta_{kj}$ describes the average rate at which infected individuals in group $j$ infect individuals of group $k$ if group $k$ were fully susceptible, $\gamma$ is the rate at which individuals recover, and $t_k$ is the rate at which we test individuals in group $k$.

In \eqref{eq:SIR-1}, a group-$k$ infected person infects people in group $k$ at rate $\beta_{kj} S_k / N_k$, since only a fraction $S_k/N_k$ of people in group $k$ are susceptible. Multiplying by the number of group-$j$ infected individuals, $I_j$, and summing across $j$ gives the rate at which people are removed from the susceptible compartment in \eqref{eq:SIR-1} and added to the infected compartment in \eqref{eq:SIR-2}. 
Each infected person recovers at a rate $\gamma$ and is tested at a rate $t_k$, resulting in $(\gamma+t_k) I_k$ being subtracted in 
\eqref{eq:SIR-2} and added in \eqref{eq:SIR-3}.


Our decision vector $x=(x_1,x_2,x_3)$, which is the amount of testing we do, is limited by a capacity constraint $\sum_{j=1}^3 N_k x_k \leq T$ for some fixed T.

We model $\beta_{kj} = \kappa\, \theta_{kj}$, where $\kappa$ is a global infectivity parameter and $\theta_{kj}$ is the number of individuals of group $k$ that an  infected individual of group $j$ interacts with per day.

The uncertain parameter vector is $\theta = (\theta_{kj} : k, j)$, which can be estimated before the university or city reopens from contact tracing data that can be obtained, with effort, from another  university or city with similar social groups. 
Contact tracing an infected person in social group $j$ yields a sample of the vector $y_{\cdot,j} = (y_{kj} : k)$, i.e. the numbers of people from each group with which that individual enjoyed social contact. We assume that $y_{kj} \sim \text{Lognormal}(\log(\theta_{kj}) - \frac{1}{2}, 1)$, and each $y_{kj}$ is independent. Letting $M_j$ be the number of contact traces of infected individuals in group $j$, our estimate for $\theta_{kj}$ is the sample mean of the observed samples for that element, i.e. $\hat{\theta}_{kj} =  \frac{\sum_{i=1}^{M_j} y_{kj}^i}{M_j}$, where $y_{\cdot,j}^i$ is the $i$th trace. Contact tracing requires time from a trained individual and so is limited. 
We write this constraint as $\sum_{j=1}^3 M_j \leq C$.

Our objective function $f(x,\theta)$ is the cumulative number of infected individuals over a finite horizon of 100 days, as given by the output of a simulation. We solve this for any fixed $\theta$ by optimizing over $x$ using L-BFGS-B. Partial derivatives in $D$ are calculated using finite differences. The optimized allocation is found via a closed form similar to Equation \ref{eq:closed}, and 1000 samples are used for optimization of the experiment design and the evaluation of the regret. The baseline allocates equal contact traces to the three groups. For more implementation details, refer to the supplement.

Results are shown in Table \ref{tab:res_pan} for parameters 
\[ \theta= \begin{pmatrix} 12 & 10 & 1\\ 10 & 8 & 1 \\ 1 & 1 & 1 \end{pmatrix},\ \kappa = \frac{1}{105},\ \gamma = 0.1,\ T=100 \]
and varying $C$, where an infected person arrives in group 2 at $T=0$. 
We use an independent Gamma prior with shape $\theta_{jk}$ and scale 1 for each $j,k$, The mean for the prior was chosen so that the first two groups have significantly more contacts than the third, modeling heterogeneity in contact rates across social groups \citep{mossong2008social}.

Higher contact rates make it more important to provide sufficient testing in the first two groups compared to the third, and also make it more important to estimate the contact rates well in these two groups. 
Indeed, the optimized allocation 
allocates more contact traces to the subpopulations with more contacts and has significantly less regret.

\begin{table}\centering
  \begin{tabular}{cccc}\toprule
  $C$ &Allocation& Opt. Regret & Uni. Regret\\
  $10$ &[5,4,1]& $37 \pm 5$ & $49 \pm 6$ \\
  $30$ &[17,12,1]& $16 \pm 2$ & $24 \pm 3$\\
  $100$ & [57,42,1]& $6.5 \pm 1.0 $ & $9.3 \pm 1.3$ \\
    $300$ & [172,126,2]& $2.2 \pm 0.3 $ & $3.1 \pm 0.4$ \\
  \bottomrule
  \end{tabular}
  \caption{Comparison of regret under optimized and baseline uniform allocations in the pandemic example of Section \ref{sec:pandemic}. The regrets reported are 95\% confidence intervals. The column ``Allocation'' shows the allocation of contact traces to groups optimizing the bound in \eqref{eq:bayes}.
  \label{tab:res_pan}}
  \end{table}
  
    \begin{figure}
  \centering
  \includegraphics[width = 0.7\textwidth]{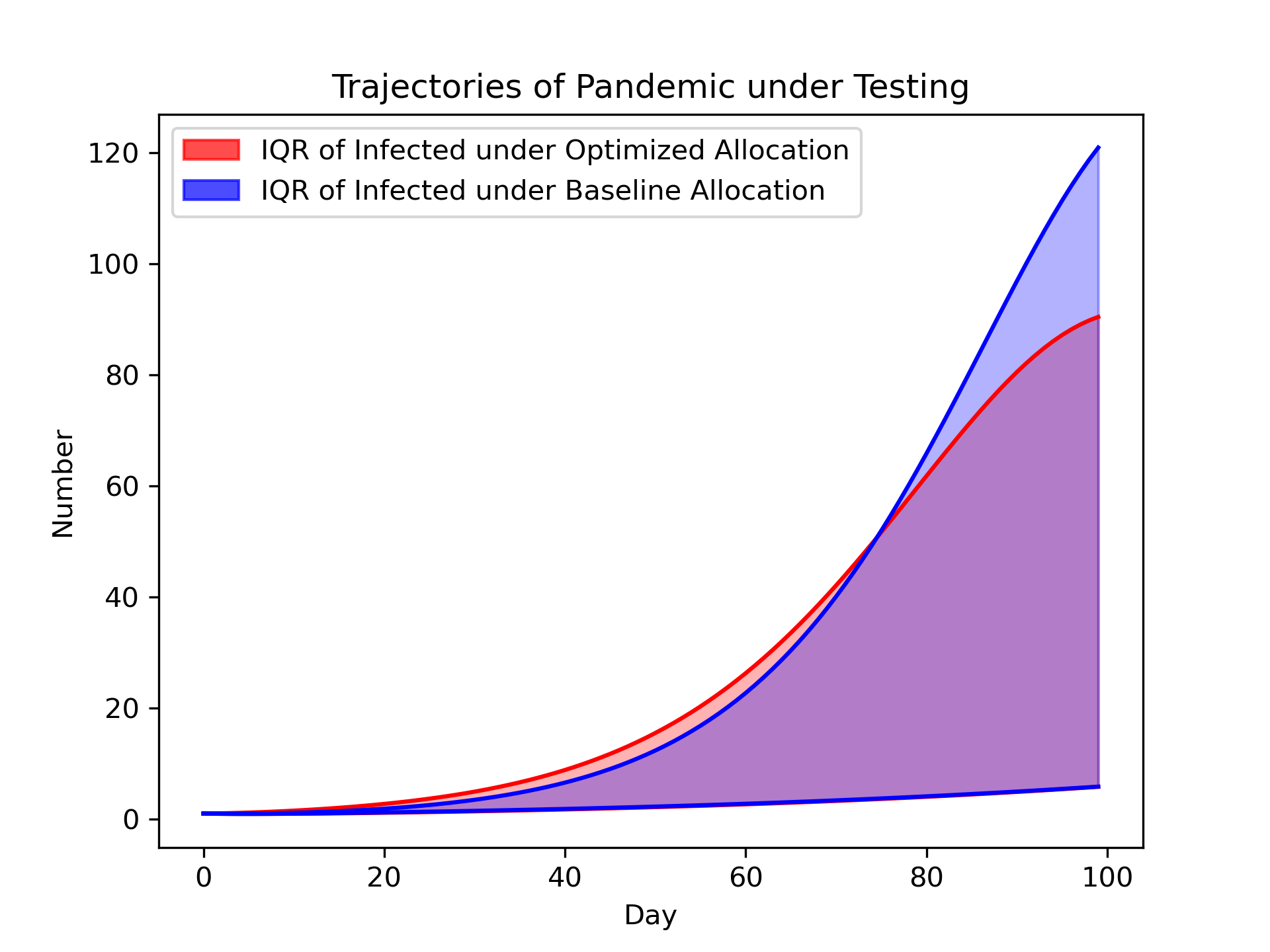}
  \caption{Interquartile range of total infected under testing decisions made based on the optimized/baseline allocation of contact traces for $C = 10$.}
  \label{fig:traj}
\end{figure}
  
  Figure \ref{fig:traj} depicts the pandemic trajectories resulting from testing decisions made under the optimized and baseline allocations for a particular sample of $\theta$. The trajectory under the optimized allocations yields notably fewer infections, demonstrating an improvement in health outcomes.
  
\begin{figure}
  \centering
  \includegraphics[width = 0.7\textwidth]{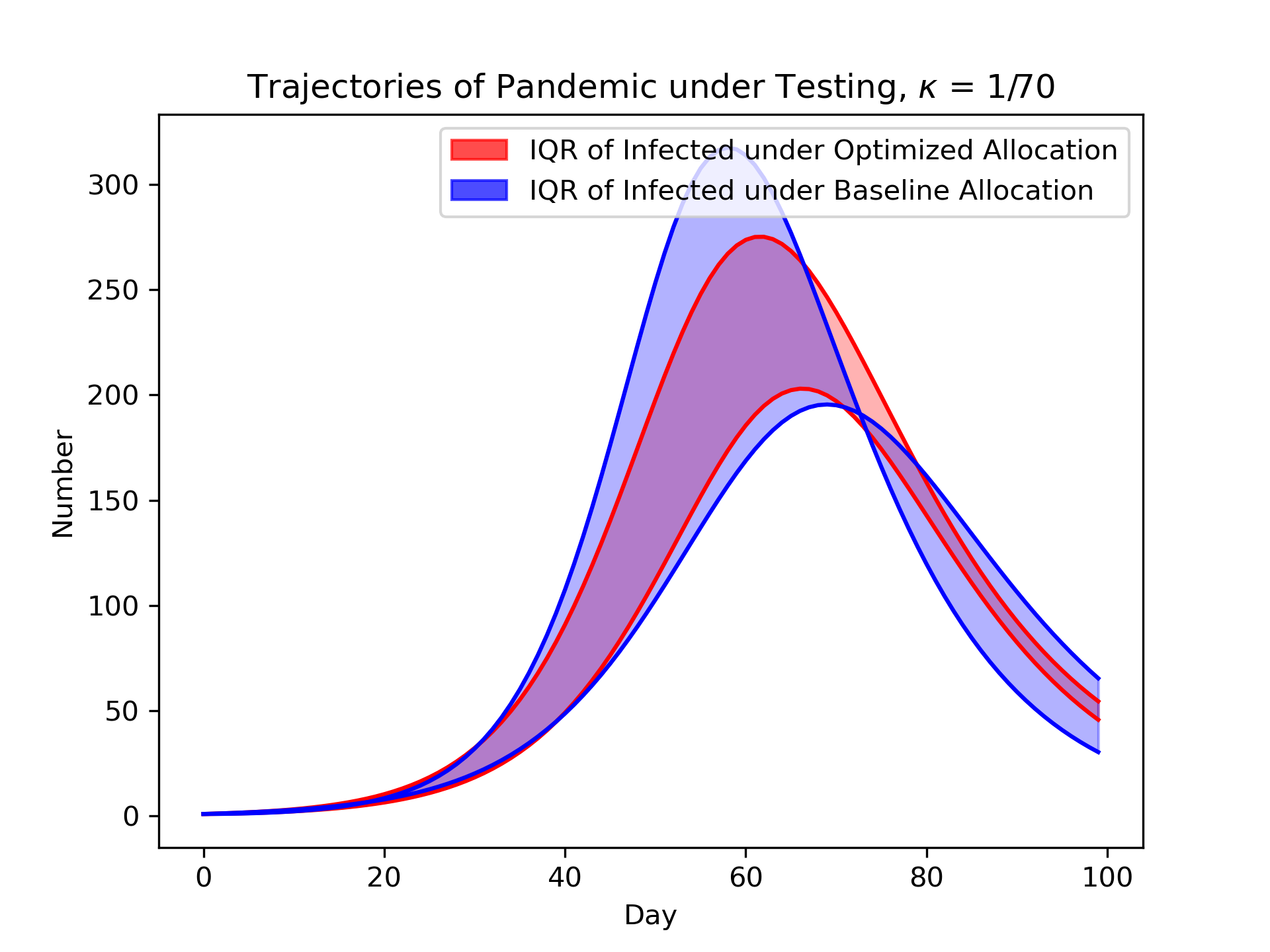}
  \caption{Interquartile range of total infected under testing decisions made based on the optimized/baseline allocation of contact traces for $C = 10$ and $\kappa = \frac{1}{70}$.}
  \label{fig:traj50}
\end{figure}

We also show results when we increase the transmissibility per contact parameter $\kappa$ by 50\%, from $\frac{1}{105}$ to $\frac{1}{70}$ (the case of decreasing the transmissibility significantly was not interesting as in either case the number of infected does go over one person). Table \ref{tab:reg50} shows the regrets for our optimized allocation versus the uniform allocation, where we change $C$, the constraint on the number of contact traces performed. Figure \ref{fig:traj50} shows the corresponding infection trajectory plots over 1000 draws from the prior. Both demonstrate that using our method results in better testing decisions that improve health outcomes, even when a larger-scale outbreak occurs.

\begin{table}\centering
  \begin{tabular}{cccc}\toprule
  $C$ &Allocation& Opt. Regret & Uni. Regret\\
  $10$ &[5,4,1]& $45 \pm 4$ & $56 \pm 5$ \\
  $30$ &[16,13,1]& $22 \pm 3$ & $26 \pm 3$\\
  $100$ & [57,42,1]& $7.6 \pm 1.1 $ & $10.1 \pm 1.2$ \\
    $300$ & [172,126,2]& $2.6 \pm 0.3 $ & $3.4 \pm 0.4$ \\
  \bottomrule
  \end{tabular}
  \caption{Comparison of regret under optimized and baseline uniform allocations in the pandemic example, when transmissibility is increased 50\%. The regrets reported are 95\% confidence intervals.}
  \label{tab:reg50}
  \end{table}
  
\section{Conclusion}
\label{sec:future}

This work introduced a novel regret bound when using estimate-then-optimize and a method for using this bound to optimize the experimental design under distributional assumptions on the underlying model parameter. To justify our theory, we analyzed three examples in detail, including a realistic application to pandemic control.

There are several directions in theory and application which would address some limitations of our method.

From a theoretical perspective, it would be interesting to consider how our approach could be extended to make decisions about whether to use only a subset of the uncertain parameters $\theta$. For instance, one may be able to conclude that because the regret is insensitive to some components of $\theta$, it may be okay to drop them completely in the experiment while still maintaining a good decision.

Similarly, one could feasibly use a similar approach for the more general problem of model selection. In this case, the $\theta$ may be infinite-dimensional in order to capture non-parametric models. As such, we imagine more technical machinery would likely be necessary.

Finally, we are improving the ability to make decisions based on data. When this is used by those whose goals are aligned with the broader society, then this should tend to improve societal outcomes, but our results could also be used by those whose goals are not so-aligned.

\bibliographystyle{unsrtnat}
\bibliography{references}  






\begin{appendix}
\section{Proofs}
\label{sec:proof}

\subsection{Proof of Theorem 3.4}
\begin{proof}
By taking the second order Taylor expansion on $f$ at $x^*, \hat{\theta}$, there exists $\overline{x} \in \mathcal{S}(\hat{x}, x^*)$ such that
\[\frac{1}{2} (\hat{x} - x^*)^\top \frac{\partial^2 f}{\partial x ^2} \Bigg|_{x= \overline{x}, \theta = \hat{\theta}} (\hat{x} - x^*)  = f(\hat{x}, \hat{\theta}) - f(x^*, \hat{\theta}) -\frac{\partial f}{\partial x} \Bigg|_{x= x^*, \theta = \hat{\theta}}^\top (\hat{x} - x^*) .\]

Observe that $\overline{x} \in B_{\hat{\theta}}$ because $B_{\hat{\theta}}$ is a convex set containing both $\hat{x}$ and $x^*$
(it contains $\hat{x} = x^*(\hat{\theta})$ by construction and it contains $x^*$ by the assumption made in the statement of the theorem), and 
$\overline{x}$ lies on a line between $x^*$ and $\hat{x}$. 

By strong convexity, and because 
$\overline{x} \in B_{\hat{\theta}}$, we have that
\begin{align*} 
\frac{1}{2} (\hat{x} - x^*)^\top \frac{\partial^2 f}{\partial x ^2} \Bigg|_{x= \overline{x}, \theta = \hat{\theta}} (\hat{x} - x^*) &\geq \rho  \| \hat{x} - x^* \|^2\\
\end{align*}

Combining the previous two equations, we have
\begin{align}
\frac{\rho }{2}  \| \hat{x} - x^*\|^2 &\leq f(\hat{x},\hat{\theta}) - f(x^*, \hat{\theta})- \frac{\partial f}{\partial x} \Bigg|_{x= x^*, \theta = \hat{\theta}}^\top (\hat{x} - x^*)\nonumber \\
&\leq f(\hat{x},\hat{\theta}) - \min_{x'} f(x', \hat{\theta}) - \frac{\partial f}{\partial x} \Bigg|_{x= x^*, \theta = \hat{\theta}}^\top (\hat{x} - x^*)\nonumber \\
&=  -\frac{\partial f}{\partial x} \Bigg|_{x= x^*, \theta = \hat{\theta}}^\top (\hat{x} - x^*)
\end{align}

Now consider the second order Taylor expansion on $\frac{\partial f}{\partial x}$ with respect to $\theta$, wherein we have that there
exists $\overline{\theta} \in \mathcal{S}(\hat{\theta}, \theta^*)$ such that 

\begin{align*}
	&-\frac{\partial f}{\partial x} \Bigg|_{x= x^*, \theta = \hat{\theta}}^\top (\hat{x} - x^*)  \\
	& =  -\frac{\partial f}{\partial x} \Bigg|_{x= x^*, \theta = \theta^*}^\top (\hat{x} - x^*)  - (\hat{\theta} - \theta^*)^\top \frac{\partial^2 f}{\partial x \partial \theta}\Bigg|_{x = x^*, \theta = \theta^*}^\top (\hat{x} - x^*) \\
	& \quad - \frac{1}{2} D^2_\theta D_x f(x^*, \overline{\theta}) [\hat{x} - x^*, \hat{\theta} - \theta^*, \hat{\theta} - \theta^*]
\end{align*}

Using first-order optimality of $x^*$, we can express this exactly as 
\[-\frac{\partial f}{\partial x} \Bigg|_{x= x^*, \theta = \hat{\theta}}  = (\hat{\theta} - \theta^*)^\top \frac{\partial^2 f}{\partial x \partial \theta}\Bigg|_{x = x^*, \theta = \theta^*}^\top (\hat{x} - x^*) - \frac{1}{2} D^2_\theta D_x f(x^*, \overline{\theta}) [\hat{x} - x^*, \hat{\theta} - \theta^*, \hat{\theta} - \theta^*] \] 

which by the second-order smoothness assumption implies

\begin{align*} -\frac{\partial f}{\partial x} \Bigg|_{x= x^*, \theta = \hat{\theta}} (\hat{x} - x^*) &\leq (\hat{\theta} - \theta^*)^\top \frac{\partial^2 f}{\partial x \partial \theta}\Bigg|_{x = x^*, \theta = \theta^*}^\top (\hat{x} - x^*) +  \frac{\beta_2}{2} \| \hat{x} - x^*\| \|\hat{\theta} - \theta^* \|^2 \\
\end{align*}

Combining these with equation (1) yields
\[ \frac{\rho}{2} \| \hat{x} - x^* \| ^2 \leq - (\hat{\theta} - \theta^*)^\top \frac{\partial^2 f}{\partial x \partial \theta}\Bigg|_{x = x^*, \theta = \theta^*}^\top (\hat{x} - x^*) + \frac{\beta_2}{2} \|\hat{x} - x^* \| \|\hat{\theta} - \theta^* \|^2 .\]

As the left is positive, so is the right, and we can square both sides to yield a valid inequality:

\begin{align*} \frac{\rho^2}{4}  \| \hat{x} - x^* \| ^4 &\leq \left( - (\hat{\theta} - \theta^*)^\top \frac{\partial^2 f}{\partial x \partial \theta}\Bigg|_{x = x^*, \theta = \theta^*}^\top (\hat{x} - x^*) + \frac{\beta_2}{2} \|\hat{x} - x^* \| \|\hat{\theta} - \theta^* \|^2 \right)^2 \\
&\leq 2 \left((\hat{\theta} - \theta^*)^\top \frac{\partial^2 f}{\partial x \partial \theta}\Bigg|_{x = x^*, \theta = \theta^*}^\top (\hat{x} - x^*) \right)^2 + \frac{\beta^2_2}{2} \| \hat{x}-x^* \|^2 \| \hat{\theta} - \theta^* \|^4\\
&\leq 2 \|\hat{x} - x^* \|^2 \|\frac{\partial^2 f}{\partial x \partial \theta}\Bigg|_{x = x^*, \theta = \theta^*} (\hat{\theta} - \theta^*) \| ^2 +  \frac{\beta^2_2}{2} \| \hat{x}-x^* \|^2 \| \hat{\theta} - \theta^* \|^4,
\end{align*}

where the first inequality uses that $(a-b)^2 \leq 2a^2 + 2b^2$ and the second uses Cauchy-Schwarz.

As $\hat{x} = x^*$ is not meaningful to consider, we assume $\hat{x} \neq x^*$ so that dividing through by $(4/\rho^2)\| \hat{x} - x^*\|^2>0$ yields

\begin{equation} \|\hat{x} - x^*\|^2 \leq \frac{8}{\rho^2} \left( \|\frac{\partial^2 f}{\partial x \partial \theta}\Bigg|_{x = x^*, \theta = \theta^*} (\hat{\theta} - \theta^*) \| ^2 + \frac{\beta^2_2}{4} \|\hat{\theta} - \theta^* \|^4 \right).
\end{equation}

Now we use a final Taylor expansion to conclude that there exists $\overline{x} \in \mathcal{S}(\hat{x}, x^*)$ such that 

\begin{align} f(\hat{x}, \theta^*) - f(x^*, \theta^*) &= \frac{\partial f}{\partial x} \Bigg|_{x = x^*, \theta = \theta^*} (\hat{x} - x^*) + \frac{1}{2} \frac{\partial^2 f}{\partial x^2} \Bigg|_{x = \overline{x}, \theta = \theta^*} (\hat{x} - x^*)^2  \nonumber \\
	&\leq \frac{\beta_1}{2} \|\hat{x} - x^* \|^2.
\end{align}

Combining (2) with (3) yields

\begin{equation}f(\hat{x}, \theta^*) - f(x^*, \theta^*) \leq  \frac{4\beta_1}{\rho^2} \left( \|\frac{\partial^2 f}{\partial x \partial \theta}\Bigg|_{x = x^*, \theta = \theta^*} (\hat{\theta} - \theta^*) \| ^2 + \frac{\beta^2_2}{4} \|\hat{\theta} - \theta^* \|^4 \right). \label{eq:main}  \end{equation}
\end{proof}

\subsection{Proof of Lemma 4.3}
\begin{proof}
We observe that under assumption 4.1, we have that

\[\sqrt{n} \Sigma^{-1/2} (\hat{\theta} - \theta^*) \sim N(0,I_d), \]
where $\theta \in \mathbb{R}^d$ and $I_d$ is the $d$-dimensional identity matrix.

We then appeal to equation (3.7) in \cite{vershynin2018high}, which gives a high probability bound for the norm of a multivariate Gaussian.

In particular, we have that 

\begin{align*}
    \mathbb{P}\left( \| \hat{\theta} - \theta^*\|_2 \geq \frac{\sqrt{\lambda_{max} (\Sigma)}}{\sqrt{n}} (t + \sqrt{d})  \right) & \leq \mathbb{P} (\|\sqrt{n} \Sigma^{-1/2} (\hat{\theta} - \theta^*) \|_2 \geq t + \sqrt{n} ) \\ 
    &= \mathbb{P}(\| g \|_2 \geq t +\sqrt{d}) \\ 
    &\leq \mathbb{P}(|\|g\|_2 - \sqrt{d}| \geq t) \\ 
    &\leq 2\exp(-ct^2)
    \end{align*}
for a universal constant $c$ and where $g \sim N(0, I_d)$.

Under assumption 4.2, we first observe that trivially $\theta^* \in \Theta^*$. Furthermore, by definition of openness we know there exists some $\epsilon$-ball $N_\epsilon(\theta^*) \subset \Theta^*$. Since we desire that $\hat{\theta}$ lies in this ball, we solve for $t$ such that 

\[ \frac{\sqrt{\lambda_{max} (\Sigma)}}{\sqrt{n}} (t + \sqrt{d})  = \epsilon \] to yield the high probability event that $\hat{\theta}$ lying in this ball.

This yields 

\[t = \frac{\epsilon \sqrt{n}}{\sqrt{\lambda_{max}(\Sigma)}} - \sqrt{d}.\]

For $n > \frac{d \lambda_{max}(\Sigma)}{\epsilon^2}$, we have that $t$ is positive for a valid bound.

Plugging this into the probability bound, we get the probability that $\|\hat{\theta} - \theta^*\| \geq \epsilon$ is bounded above by $2\exp(-cd) \exp\left(- \frac{cn\epsilon^2}{\lambda_{max}(\Sigma)} + \frac{2c\epsilon \sqrt{nd}}{\sqrt{\lambda_{max}(\Sigma)}} \right)$.

Taking $n$ larger than, say, $\frac{16d \lambda_{max}(\Sigma)}{\epsilon^2}$, we have an upper bound of
$2\exp(-cd) \exp\left(- \frac{cn\epsilon^2}{2\lambda_{max}(\Sigma)} \right)$.

So for $n$ sufficiently large, i.e. $n \geq \max \left(\frac{d \lambda_{max}(\Sigma)}{\epsilon^2}, \frac{16d \lambda_{max}(\Sigma)}{\epsilon^2} \right)  = \frac{16d \lambda_{max}(\Sigma)}{\epsilon^2}$, taking the complementary event concludes the proof.



\end{proof}

\subsection{Proof of Theorem 4.4}
\begin{proof}
By the previous theorem,
$f(\hat{x}, \theta^*) - f(x^*, \theta^*) \le \frac{4\beta_1}{\rho^2} \left( \|D (\hat{\theta} - \theta^*) \| ^2 + \frac{\beta^2_2}{4} \|\hat{\theta} - \theta^* \|^4 \right)$
when 
$x^* \in B_{\hat{\theta}}$.

Therefore, the probability that 
\begin{align*} f(\hat{x}, \theta^*) - f(x^*, \theta^*)
> &\frac{4\beta_1}{\rho^2}  \left( \Tr \frac{D\Sigma D^\top}{n} + 2 \sqrt{\Tr \left[ \left(\frac{D\Sigma D^\top}{n}\right)^2\right] \log n } + 2 \|\frac{D\Sigma D^\top}{n}\|_2 \log n \right)\\
&+ o(\Tr \frac{D\Sigma D^\top}{n} \log n)
\end{align*} is bounded above by the sum of 
\begin{itemize}
\item the probability that $x^* \notin B_{\hat{\theta}}$.
\item the probability that $ \|D (\hat{\theta} - \theta^*) \| ^2 \geq  \Tr \frac{D\Sigma D^\top}{n} + 2 \sqrt{\Tr \left [\left(\frac{D\Sigma D^\top}{n}\right)^2 \right] \log n} + 2 \|\frac{D\Sigma D^\top}{n}\|_2 \log n $,
\item the probability that $\|\hat{\theta} - \theta^* \| ^4 \in  o(\Tr \frac{D\Sigma D^\top}{n} \log n)$.
\end{itemize}

By the lemma, the probability that 
$x^* \notin B_{\hat{\theta}}$ is bounded above by $\alpha \exp(-\beta n)$.

By \cite{hsu2012tail}, the probability that

\[ \|D (\hat{\theta} - \theta^*) \| ^2 \geq \Tr \frac{D\Sigma D^\top}{n} + 2 \sqrt{\Tr \left[ \left(\frac{D\Sigma D^\top}{n}\right)^2 \right]} \sqrt{t} + 2 \|\frac{D\Sigma D^\top}{n}\|_2 t \] 

is bounded above by $\exp(-t).$ Choosing $t = \log n$ yields the first term in bound.

For the little-$o$ term, the same bound from \cite{hsu2012tail} yields 

\[ \|\hat{\theta} - \theta^* \| ^4 \geq \left( \Tr \frac{\Sigma}{n} + 2 \sqrt{\Tr \left(\frac{\Sigma }{n}\right)^2 } \sqrt{t} + 2 \|\frac{\Sigma}{n}\|_2 t \right)^2 \]
with probability $\exp(-t)$. Taking the same $t = \log n$ yields a term which is $o(\Tr \frac{D\Sigma D^\top}{n} \log n)$ (note $\Tr D \Sigma D^\top$ and $\Tr \Sigma$ are within a constant factor corresponding to eigenvalues of $D$.)


Then by the union bound, the probability that 
$f(\hat{x}, \theta^*) - f(x^*, \theta^*)
> \Tr \frac{D\Sigma D^\top}{n} + 2 \sqrt{\Tr \left[ \left(\frac{D\Sigma D^\top}{n}\right)^2\right] }\log n  + 2 \|\frac{D\Sigma D^\top}{n}\|_2 \log n + o(\Tr \frac{D\Sigma D^\top}{n} \log n)$ 
is bounded above by $\alpha \exp(-\beta n) + \frac{2}{n}$. Taking the complementary event completes the proof.

Finally, we note that, as a function of $D\Sigma D^\top$, the first trace term is the dominant term because $D\Sigma D^\top$ is symmetric and PSD, and therefore all eigenvalues are nonnegative. Therefore the sum of the eigenvalues (the trace), or equivalently the 1-norm of the eigenvalues, is larger than the 2-norm (the second trace term), and also larger than the largest eigenvalue (the 2-norm term).

\end{proof}

\begin{remark}
The choice of $t = \log n$ in the bound above is due to the following argument. Suppose $\theta$ lies in some compact set. Then there exists some upper bound on the regret, say $M$. Then for large $n$, an upper bound on the expected regret reduces to $1 \cdot \frac{t}{n} + M \cdot \exp(-t).$ Minimizing this as a function of $t$ yields $t = \log \frac{R}{c}  + \log n$.
\end{remark}

\subsection{Verifying Assumptions for Pricing Example}

One may easily derive the following

\begin{equation}\frac{\partial f}{\partial x} = -c(x,\theta) + \theta_1 x c(x,\theta)(1-c(x,\theta)). \label{eq:price}\end{equation}

To verify the convexity assumptions, first note that setting 
\eqref{eq:price} to 0 yields that the optimal solution $x^*$ satisfies 
\[x^* = \frac{1}{\theta_1 (1 - c(x^*,\theta))}.\]
(This is not a closed form solution as the RHS depends on $x$, but it shows that $x^*>0$.)

Furthermore, the second derivative of $f$ with respect to $x$ is

\begin{equation} \frac{\partial^2 f}{\partial x^2} = c(x,\theta)(1-c(x,\theta)) \left( 1 + \theta_1 - \theta_1^2 x(1-2c(x,\theta)) \right) \label{eq:price2nd}\end{equation}

Therefore we get
\[\frac{\partial^2 f}{\partial x^2}\Bigg|_{x = x^*} = c(x^*,\theta)(1-c(x^*,\theta)) \left(1 + \theta_1 - \theta_1 \frac{1 - 2c(x^*,\theta)}{1- c(x^*,\theta)} \right) > 0, \]

since we assumed $\theta_1>0$. Also, for any $\theta$, it is easy to see that this is bounded away from 0 for $x$ in a neighborhood around $x^*$.

Furthermore, we observe that Equation \eqref{eq:price2nd}
is clearly bounded for $x \in [0,\infty)$ and for all $\theta$ since the logistic function decays exponentially.

Finally, it remains to check the Hessian of \eqref{eq:price} with respect to $\theta$.
One can straightforwardly derive the expressions for the terms of the Hessian:
\begin{align*}
    \frac{\partial f^3}{\partial x \partial \theta_0^2} =& c(x,\theta)(1-c(x,\theta))(2c(x,\theta) - 1) \\
    & + \theta_1 xc(x,\theta)(1-c(x,\theta)) \left((2c(x,\theta)-1)^2 + 2c(x,\theta)(1-c(x,\theta))\right)\\
    \frac{\partial f^3}{\partial x \partial \theta_0 \partial \theta_1 } =& 2x c(x,\theta)(1-c(x,\theta))(2c(x,\theta) - 1) \\
    & + \theta_1 x^2 c(x,\theta)(1-c(x,\theta)) \left((2c(x,\theta)-1)^2 + 2c(x,\theta)(1-c(x,\theta))\right)\\
    \frac{\partial f^3}{\partial x \partial \theta_1^2} =& 3x^2 c(x,\theta)(1-c(x,\theta))(2c(x,\theta) - 1) \\
    & + \theta_1 x^3 c(x,\theta)(1-c(x,\theta)) \left((2c(x,\theta)-1)^2 + 2c(x,\theta)(1-c(x,\theta))\right)
\end{align*}

Each term in the Hessian includes a multiplicative factor of $c(x,\theta)(1-c(x,\theta)$. Therefore for any $x\in [0,\infty)$ and any $\theta$,
each term in the Hessian, and thus also its largest eigenvalue, must remain bounded as the logistic function decays exponentially.

\section{Implementation Details for Pandemic Example}

For our pandemic application, we built a discrete-time simulation of the multi-group SIR model as described in the paper. In particular, for time $t = 1, \ldots, 100$, we have that 

\begin{align*}
    S[t] &= S[t-1] - S[t-1] \circ \beta I[t-1] \circ \frac{1}{N}\\
    I[t] &= I[t-1] + S[t-1] \circ \beta I[t-1] \circ \frac{1}{N} - (\gamma\mathbf{1} + x)\circ I[t-1] \\
    R[t] &= R[t-1] + (\gamma\mathbf{1} + x)\circ I[t-1],
\end{align*}
where $\circ$ denotes entrywise product and $\mathbf{1}$ denotes a vector of ones. We treat $f(x,\theta)$ as the output of the simulation, where recall that $\beta = \kappa \theta$.

When optimizing $f$ with respect to the decision variable $x$ for a given $\theta$, we use L-BFGS-B. Since L-BFGS-B is used for unconstrained problems, we reparameterized the problem to only have bounds on the variables instead of constraints. In particular, we define another set of variables $y_1, y_2$, where $y_1$ is the proportion of testing capacity used by the first group, $y_2$ is the proportion of the testing capacity which remains that the second group uses. Therefore we may include $y_1, y_2 \in [0,1]$, which is amenable to usage by L-BFGS-B.  

We use finite difference estimates for the derivative term $D$, with step size $h = 10^{-6}$. In particular, we use
\[\frac{\partial^2 f}{\partial x \partial \theta} \approx \frac{f(x+h, \theta+h) - f(x+h, \theta-h) - f(x-h, \theta+h) + f(x-h, \theta-h)}{4h^2}.\] 

Finally, $D$, the trace term is minimized in closed-form. In particular, it is easy to see that we have a problem of the form (after relaxing integrality)
\begin{align*}
    \min_{n_1, \ldots, n_3} \, & \sum_{i=1}^3 \frac{\rho_i^2}{M_i} \\ 
    \text{s.t.} & \sum_{i=1}^3 M_i = C \\ 
    & M_i \geq 0 \quad i = 1, 2, 3
\end{align*}
where $\rho_i^2 = \sum_{\ell=1}^2 \sum_{j=1}^3 \sum_{k=1}^3 \frac{\partial^2 f}{\partial x_\ell \partial \theta_{kj}} \sigma_{kj}^2$ for $\sigma_{kj}^2$ the variance of a $\text{Lognormal}(\log(\theta_{kj}, 1))$ random variable.

The Karush-Kuhn-Tucker conditions yield the closed-form solution

\[M_i = \frac{C \rho_i}{\sum_{j=1}^d \rho_j}, \]
and we round the solution to the nearest integer (rounding the first two and taking the difference from $C$ as $M_3$). If $M_3$ is 0 (which we disallow as then no samples are collected), we set $M_3 = 1$ and decrement one of $M_1, M_2$ each with probability 0.5.

\section{Additional Experiments}
\label{sec:add}

\subsection{Pricing}

\begin{figure}
      \centering
      \includegraphics[width = 0.7\textwidth]{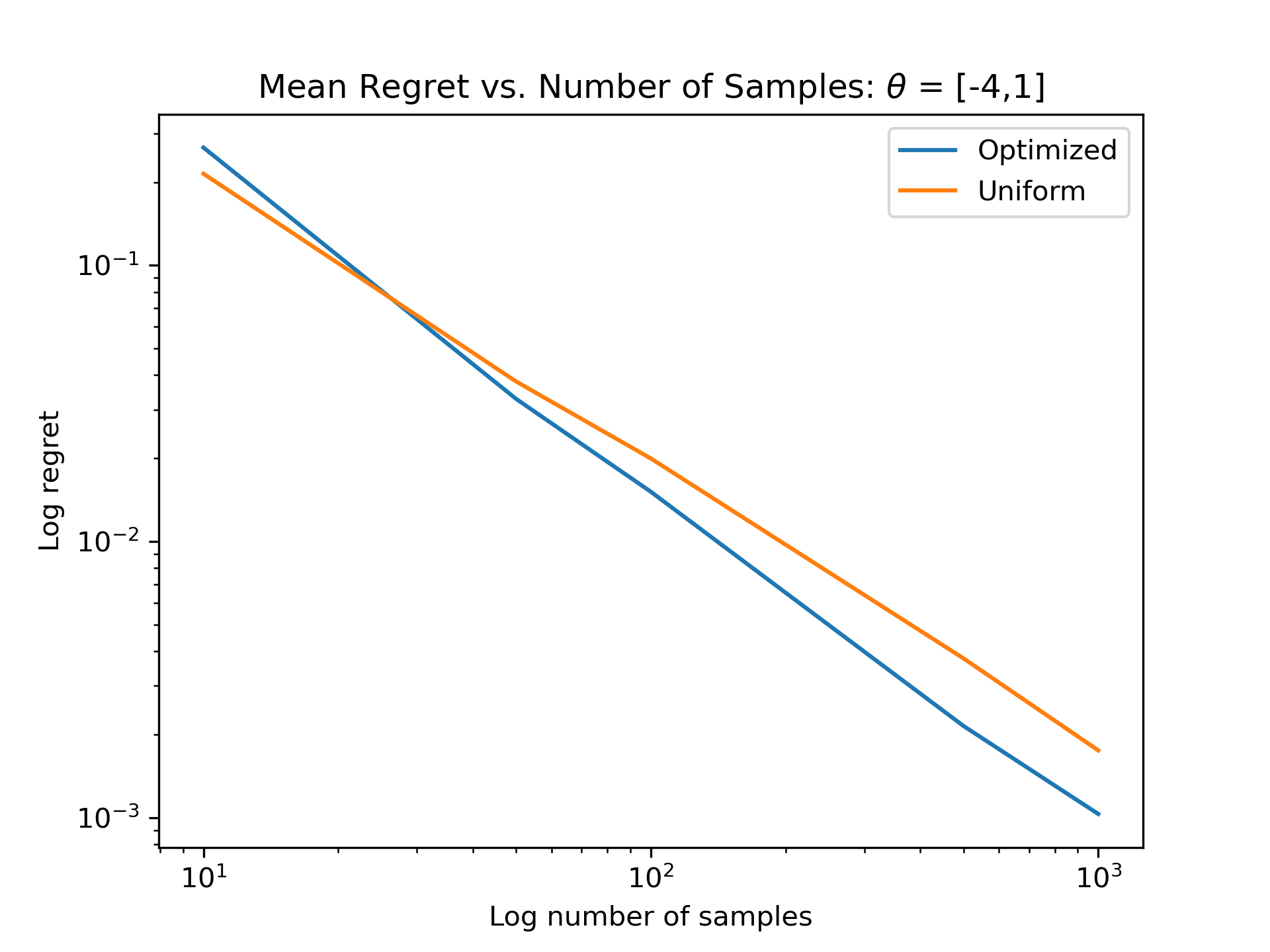}
      \caption{Regret in the pricing example decays as $O(\frac{1}{n})$, where $n$ is the number of samples allocated. Here $\theta = [-4,1]$.}
      \label{fig:n}
  \end{figure}

Figure \ref{fig:n} shows regret versus the number of samples used for the allocation, for a particular choice of $\theta = [-4,1]$. The linear relationship on the log-log scale confirms that the regret behaves as $\frac{1}{n}$. This aligns with our regret bound behaving as $O(\frac{\log n}{n})$, where the $\log n$ is obscured.

 We also see from the above graph that the relative improvement of our optimal allocation improves as $n$ increase. Table \ref{tab:res1k} shows the regret of the optimal and uniform allocations for $n=1000$ samples and more values of $\theta$, showing further relative improvement relative to when $n=100$.
 
 \begin{table}\centering
  \begin{tabular}{cccc}\toprule
  &Optimized Allocation & Uniform Allocation \\
  $\theta = [-1,1]$ & $1.0 \pm 0.3$ & $1.6 \pm 0.4$ \\
  $\theta = [-2,1]$ & $1.0 \pm 0.3$ & $2.0 \pm 0.5$\\
  $\theta = [-3,1]$ & $0.9 \pm 0.2 $ & $1.8 \pm 0.4$ \\
  $\theta = [-4,1]$ & $1.3 \pm 0.4 $ & $1.7 \pm 0.4$ \\
  $\theta = [-5,1]$ & $1.1 \pm 0.3 $ & $2.4 \pm 0.6$ \\
  $\theta = [-6,1]$ & $1.7 \pm 0.5 $ & $1.9 \pm 0.4$ \\
  $\theta = [-7,1]$ & $1.6 \pm 0.4$ & $2.6 \pm 0.7$ \\
  $\theta = [-8,1]$ & $3.9 \pm 1.5$ & $6.7 \pm 4.4$\\
  \bottomrule
  \end{tabular}
  \caption{Comparison of regret under optimized and baseline uniform allocations in the pricing example, where 1000 samples are allocated. The quantities are 95\% confidence intervals, where the units are in $10^{-3}$. 
  \label{tab:res1k}}
  \end{table}

\subsection{Pandemic}

\begin{figure}
      \centering
      \includegraphics[width = 0.7\textwidth]{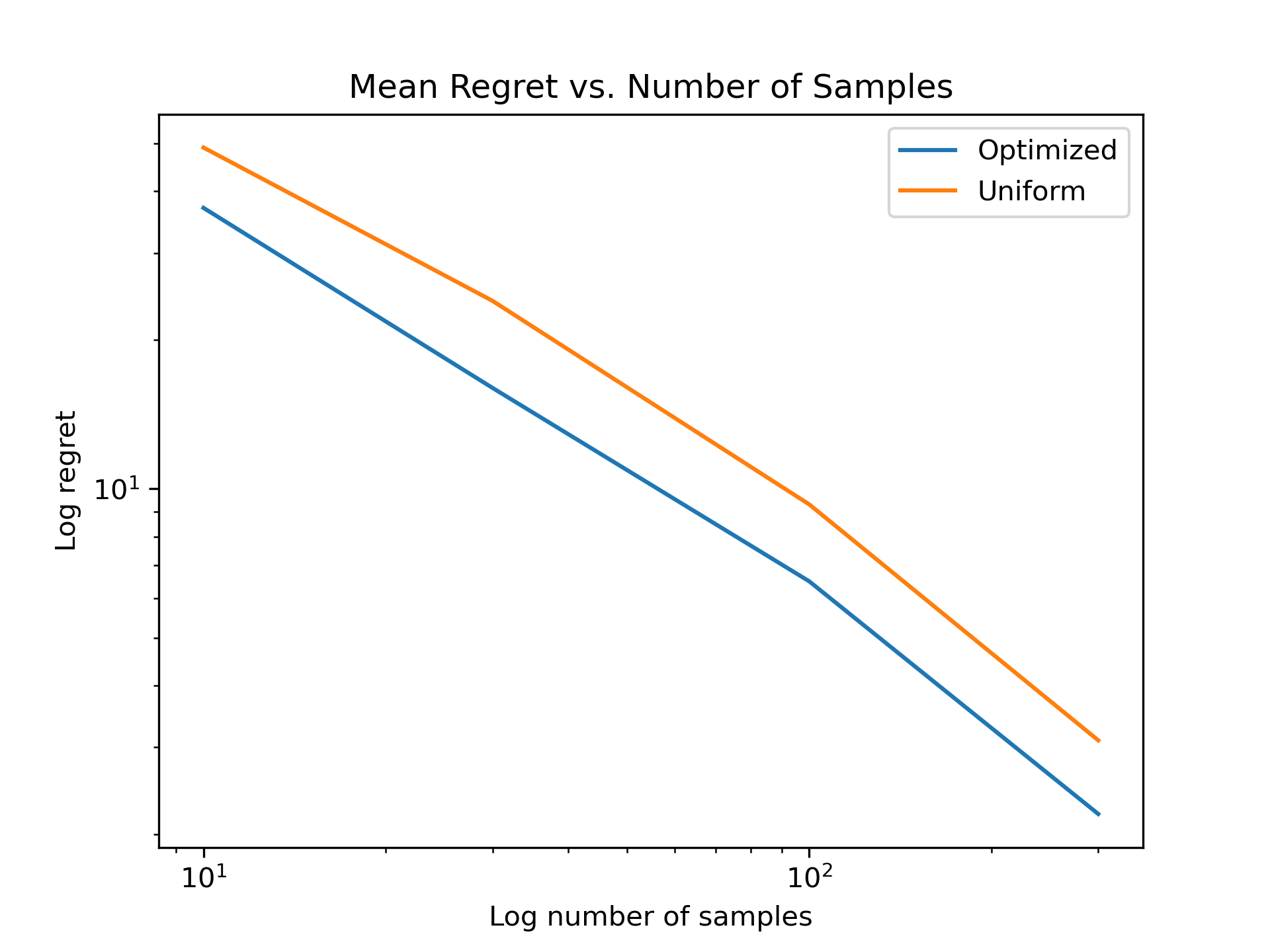}
      \caption{Regret in the pandemic example also roughly decays as $O(\frac{1}{n})$, where $n$ is the number of samples allocated.}
      \label{fig:C}
  \end{figure}

Figure \ref{fig:C} shows the analogous plot of regret vs number of samples allocated to Figure \ref{fig:n} from the pricing example. Again, we see the roughly linear relationship on the log-log scale confirms that the regret behaves as $\frac{1}{n}$, agreeing with our bound.

 \section{Code}

Please find all code used at this anonymized repository:  \url{https://anonymous.4open.science/r/ExpDesign4ETO-DEB3/}.

\end{appendix}

\end{document}